\documentclass[12 pt]{article}
\usepackage[english]{babel}
\usepackage{amsmath}
\usepackage{amssymb}

\newcommand{\norm}[1]{\left\Vert #1\right\Vert}
\newcommand{\abs}[1]{\left\vert #1\right\vert}
\newcommand{\R}{   {\ifmmode{{\mathbb R}}\else{$\mathbb R$}\fi}}
\newcommand{\N}{   {\ifmmode{{\mathbb N}}\else{$\mathbb N$}\fi}}
\newcommand{\V}{\mathcal{V}}
\newcommand{\K}{\mathcal{K}}
\newcommand{\demo}{\noindent\textit{Proof. }}

\newcommand{\rank}{\mbox{rank}\ }

\newcommand{\dt}{\frac{d}{dt}}

\begin{document}

\title{\bf Geometry of the Limit Sets of Linear Switched Systems}
\author{Moussa BALDE\footnote{LGA-LMDAN Dept. de Math\'ematiques et Informatique,UCAD, Dakar-Fann, Senegal, E-mail: moussa.balde.math@ucad.edu.sn}, Philippe JOUAN\footnote{LMRS, CNRS UMR 6085, Universit\'e
    de Rouen, avenue de l'universit\'e BP 12, 76801
    Saint-Etienne-du-Rouvray France. E-mail: Philippe.Jouan@univ-rouen.fr}}

\date{\today}

\maketitle

\begin{abstract}
The paper is concerned with asymptotic stability properties of linear switched systems. Under the hypothesis that all the subsystems share a non strict quadratic Lyapunov function, we provide a large class of switching signals for which a large class of switched systems are asymptotically stable. For this purpose we define what we call non chaotic inputs, which generalize the different notions of inputs with dwell time.

Next we turn our attention to the behaviour for possibly chaotic inputs. To finish we give a sufficient condition for a system composed of a pair of Hurwitz matrices to be asymptotically stable for all inputs.

\vskip 0.2cm


\vskip 0.2cm


\end{abstract}

\section{Introduction.}

Let $\{B_1,\dots,B_p\}$ be a finite collection of $d\times d$ matrices assumed to share a common quadratic Lyapunov function. In case where this Lyapunov function is strict, all the $B_i$'s are Hurwitz and the switched system
\begin{equation}
\label{switchsys}
\dot{x}=B_{u(t)}x
\end{equation}
is asymptotically stable for any switching input $u$ (see for instance \cite{liberzon-book}, and \cite{agr}, \cite{bbm}, \cite{AL00}, \cite{branicky98} for other approaches). 

In this paper we investigate the case where the Lyapunov function is not strict. For each matrix $B_i$ the space admits an orthogonal and $B_i$-invariant decomposition $\V_i\oplus \V_i^{\bot}$, such that the restriction of $B_i$ to $\V_i$ is skew-symmetric in a suitable basis and its restriction to $\V_i^{\bot}$ is Hurwitz (see Lemma \ref{decomposition} in Section \ref{Preliminaries}). In the very interesting work \cite{SVR} each matrix $B_i$ is assumed to vanish on the subspace $\V_i$.

Our aim is to obtain asymptotic stability results without this assumption. For this purpose we introduce two fundamental tools. The first one consists in lifting the problem to the space of matrices $\mathcal{M}(d;\R)$, as it is often done for the ordinary linear differential equations, and in applying the polar decomposition to the matrix trajectory. This enables us to consider a symmetric matrix $S_u$, defined as a limit (see Section \ref{lift}), that depends on the input, and is equal to zero if and only if the switched system is asymptotically stable for that last (Theorem \ref{squareroot}). On the other hand we use Ascoli's Theorem to show that the matrix trajectory converges uniformly on some sequences of intervals (Section \ref{interval}).

Thanks to these tools we define in Section \ref{regss} what we call "non chaotic inputs". They generalize the various notions of inputs with dwell-time, or with average dwell-time that can be found in the literature (see for instance \cite{Morse93},\cite{HM}, \cite{SVR}, or \cite{liberzon-book} for a general reference).

The main results are in Sections \ref{regss}, \ref{fastinputs}, and \ref{pairs}.

In Section \ref{regss} we state asymptotic stability criterions for regular inputs, that is for non chaotic and in some sense well-distributed inputs (Theorems \ref{main} and \ref{corollaire}).

Section \ref{fastinputs} deals with general, that is possibly chaotic, inputs (Theorems \ref{tjrsvrai} and \ref{fastswitch}).

A general result of asymptotic stability for pairs of Hurwitz matrices is established in Section \ref{pairs}.

Section \ref{ExandAp} is devoted to examples.


\section{The systems under consideration}\label{Preliminaries}
\subsection{The matrices}

As explained in the introduction we deal with a finite collection of $d\times d$ matrices, $\{B_1,B_2,\dots,B_p\}$, assumed to share a common, but not strict in general, quadratic Lyapunov function. More accurately there exists a symmetric positive definite matrix $P$ such that the symmetric matrices $B_i^TP+PB_i$ are nonpositive ($B^T$ stands for the transpose of $B$).
Since the Lyapunov matrix $P$ is common to the $B_i$'s we can assume without loss of generality that $P$ is the identity matrix, in other words that $B_i^T+B_i$ is non positive for $i=1,...,p\ $:
\begin{equation}\label{condition}
\forall x\in \R^d, \quad \forall i=1,...,p \qquad x^T(B_i^T+B_i)x\leq 0.
\end{equation}

\noindent{\bf Norms}. The natural scalar product of $\R^d$ in this context is the canonical one, defined by $<x,y>=x^Ty$ (it would be $x^TPy$ if the Lyapunov matrix were $P$). The norm of $\R^d$ is consequently chosen to be $\norm{x}=\sqrt{x^Tx}$ and the space $\mathcal{M}(d,\R)$ of square matrices (or equivalently of endomorphisms of $\R^d$) is endowed with the related operator norm:
$$
\norm{B}=\max\{\norm{Bx};\ x\in\R^d \mbox{ and } \norm{x}=1\}.
$$

The following lemma is known, but fundamental. For this reason we give here a proof, which in our view enlightens the sequel.

\newtheorem{decomposition}{Lemma}
\begin{decomposition} \label{decomposition}
Let $B$ be a $d\times d$ matrix, identified with an endomorphism of $\R^d$, and assumed to satisfy $B^T+B\leq 0$. Let
$$
\V=\{x\in\R^d;\ \forall t\in \R\ \norm{e^{tB}x}=\norm{x}\}
$$
Then $\V$ is a $B$-invariant subspace of $\R^d$ and the restriction of $B$ to $\V$ is skew-symmetric in any orthonormal basis. The orthogonal complement $\V^{\bot}$ of $\V$ in $\R^d$ is also $B$-invariant and the restriction of $B$ to $\V^{\bot}$ is Hurwitz.

Moreover a point $x\in\R^d$ belongs to $\V$ if and only if there exists $\tau>0$ such that $\displaystyle \norm{e^{\tau B}x}=\norm{x}$.

\end{decomposition}

\demo
First of all notice that the condition $B^T+B\leq 0$ implies $\forall x\in\R^n$, $\forall t\geq 0$, $\norm{e^{tB}x}\leq\norm{x}$, so that $B$ is Hurwitz if and only if $\V=\{0\}$.

Let $B=S+A$ be the decomposition of $B$ into a symmetric part $S$ and a skew-symmetric one $A$. We have
$$
\begin{array}{ll}
x\in\V & \Longleftrightarrow \forall t\in\R \qquad \norm{e^{tB}x}^2=\norm{x}^2\\
         & \Longleftrightarrow \forall t\in\R \qquad \dt \norm{e^{tB}x}^2=0\\
         & \Longleftrightarrow \forall t\in\R \qquad x^Te^{tB^T}(B^T+B)e^{tB}x=0\\
         & \Longleftrightarrow \forall t\in\R \qquad e^{tB}x\in \ker(S)
\end{array}
$$
because $B^T+B=2S$. This shows that $\V$ is a subspace of $\R^d$, included in $\ker(S)$, and moreover $B$-invariant, as shown by a second derivation. Let $y\in \V^{\bot}$. For all $x\in \V$, it follows from $Sx=0$ that
$$
\begin{array}{ll} 
<x,By> & = x^TSy+x^TAy=(Sx)^Ty-(Ax)^Ty\\
       & = -(Ax)^Ty= -(Sx)^Ty-(Ax)^Ty =-<Bx,y>=0.  
\end{array}
$$  
The subspace $\V^{\bot}$ is therefore $B$-invariant, the restriction of $B$ to $\V^{\bot}$ is Hurwitz (if not the intersection with $\V$ would not be $\{0\}$), and its restriction to $\V$ is skew-symmetric in any orthonormal basis because $\V$ is included in the kernel of $S$.

Assume to finish that $\displaystyle \norm{e^{\tau B}x}=\norm{x}$ for some $\tau>0$. Then the equality holds for $t\in [0,\tau]$, because  $\norm{x}=\norm{e^{\tau B}x}\leq\norm{e^{tB}x}\leq\norm{x}$ for such a $t$, and by analycity for all $t\in \R$.

\hfill $\Box$

In the sequel, the set $\{x\in\R^d;\ \forall t\in \R\ \norm{e^{tB_i}x}=\norm{x}\}$ is denoted by $\V_i$ for $i=1,\dots,p$.


\subsection{The switching signal}

An input, or switching signal, is a piecewise constant and right-continuous function $u$ from $[0,+\infty[$ into $\{1,...,p\}$. We denote by $(a_n)_{n\geq 0}$ the sequence of switching times (of course $a_0=0$ and the sequence is strictly increasing to $+\infty$). Therefore $u(t)$ is constant on each interval $[a_n,a_{n+1}[$, and $u_n\in \{1,...,p\}$ will stand for this value. The duration $a_{n+1}-a_n$ is denoted by $\delta_n$.

As the input is entirely defined by the switching times and the values taken at these instants we can write
$$
u=(a_n,u_n)_{n\geq 0}.
$$

Such a switching signal being given the switched system under consideration is the dynamical system defined in $\R^d$ by
\begin{equation}\label{system}
\dot{x}=B_{u(t)}x
\end{equation}
Its solution is, for the initial condition $x$ and for $t\geq 0$,
$$
t\longmapsto \Phi_u(t)x
$$
where
$$
\Phi_u(t)=\exp((t-a_n)B_{u_n})\exp(\delta_{n-1}B_{u_{n-1}})\dots\exp(\delta_1B_{u_1})\exp(\delta_0 B_{u_0})
$$
if $t\in [a_n,a_{n+1}[$.


\subsection{The $\omega$-limit sets}
For $x\in \R^d$ we denote by $\Omega_u(x)$ the set of $\omega$-limit points of $\{\Phi_u(t)x;\ t\geq 0\}$, that is the set of limits of sequences $(\Phi_u(t_k)x)_{k\geq 0}$, where $(t_k)_{k\geq 0}$ is strictly increasing to $+\infty$.

Thanks to Condition (\ref{condition}), the norm $\norm{\Phi_u(t)x}$ is nonincreasing, and from this fact we can easily deduce the proposition:

\newtheorem{P1}{Proposition}
\begin{P1}\label{P1}
For any initial condition $x$ the $\omega$-limit set $\Omega(x)$ is a compact and connected subset of a sphere $S_r=\{x\in \R^d;\ \norm{x}=r\}$ for some $r\geq 0$.
\end{P1}

\demo The trajectory $\{\Phi_u(t)x;\ t\geq 0\}$ is bounded by $\norm{x}$. It is a general fact that the $\omega$-limit set of a bounded trajectory in a finite dimensional space is compact and connected (see for instance \cite{Marle}). Let us prove that this set is contained in a sphere. Pick two limit points,
$$
l=\lim_{k\longmapsto +\infty}\Phi_u(t_k)x \quad\mbox{ and }\quad l'=\lim_{j\longmapsto +\infty}\Phi_u(t_j)x.
$$
For any $k\geq 0$ there exists $j$ such that $t_j\geq t_k$. But $\norm{\Phi_u(t)x}$ being nonincreasing this implies $\norm{\Phi_u(t_j)x}\leq \norm{\Phi_u(t_k)x}$. We have therefore $l'\leq l$, and the converse as well.

\hfill $\Box$

\noindent{\bf Remark}. Proposition \ref{P1} is actually proved in \cite{SVR} but under an additional assumption of "paracontraction".


\section{The lift of the problem}\label{lift}

Let us now forget the initial condition $x\in\R^d$ and deal with the matrix function $\Phi_u$. That last is continuous on $[0,+\infty[$ and takes its values in the closed ball $K=B'(0,1)=\{M\in \mathcal{M}(d;\R);\ \norm{M}\leq 1\}$ of $\mathcal{M}(d;\R)$.

Its $\omega$-limit set, that we denote by $\Omega_u$, is therefore a compact and connected subset of $K$.

Let $x\in \R^d$.
If $M\in \Omega_u$ it is clear that $Mx\in \Omega_u(x)$. Conversely let $l\in \Omega_u(x)$. There exists a increasing sequence $(t_k)_{k\geq 0}$ such that $\displaystyle l=\lim_{k\longmapsto +\infty} \Phi_u(t_k)x$. For all $k$, the matrix $\Phi_u(t_k)$ belongs to the compact $K$, and we can extract a subsequence $(\Phi_u(t_{k_j}))_{j\geq 0}$ that converges to a limit $M$. We have clearly $l=Mx$, and the following proposition holds:

\newtheorem{limitset}[P1]{Proposition}
\begin{limitset}\label{limitset}
The set $\Omega_u$ is a compact and connected subset of \\$K=B'(0,1)\subset \mathcal{M}(d;\R)$, and for all $x\in\R^d$ the set $\Omega_u(x)$ is equal to:
$$
\Omega_u(x)=\{Mx;\ \ M\in \Omega_u\}=\Omega_ux.
$$
\end{limitset}

Let us now describe $\Omega_u$ more accurately. We know by Proposition \ref{P1} that the set $\Omega_ux$ is, for each $x\in \R^d$, included in a sphere:
$$
\exists r\geq 0 \quad\mbox{ such that }\ \forall M\in \Omega_u \quad \norm{Mx}=r.
$$
Pick two matrices $M$ and $N$ in $\Omega_u$. Then
\begin{equation}\label{sym}
\forall x\in \R^d \qquad x^TM^TMx=x^TN^TNx
\end{equation}
The matrix $M^TM-N^TN$ being symmetric this equality implies $M^TM-N^TN=0$, and we get
$$
\forall M,N \in \Omega_u \qquad M^TM=N^TN.
$$
Consider the polar decomposition of $M\in \Omega$: there exist an orthogonal matrix $O$ and a symmetric nonnegative one $S$ such that $M=OS$. Notice that whenever $M$ is not invertible the matrix $S$ is not definite and the matrix $O$ is not unique. However $S$ is well defined because we have
$$
M^TM=S^TO^TOS=S^2
$$
and $S$ is therefore the unique nonnegative, symmetric square root of $M^TM$ and does not depend on a particular choice of the matrix $M\in \Omega_u$. This matrix $S$ is from now on denoted by $S_u$.

We can also apply the polar decomposition to $\Phi_u(t)$ and write
$$
\Phi_u(t)=O(t)S(t) \qquad \forall t\geq 0
$$
where $O(t)$ is orthogonal and $S(t)$ symmetric positive definite. The matrix $\Phi_u(t)$ being a product of exponentials, its determinant is always positive, and $O(t)$ actually belongs to $SO_d$.
On the other hand $\Phi_u(t)^T\Phi_u(t)=S^2(t)$ decreases to $S_u^2$ as $t$ goes to $+\infty$. Indeed it is a symmetric positive matrix and for all $x\in\R^d$ the norm $x^T\Phi_u(t)^T\Phi_u(t)x=\norm{\Phi_u(t)x}^2$ is nonincreasing.

The convergence of $S(t)$ to $S_u$ has the following consequence.

Let $(\Phi_u(t_k)=O(t_k)S(t_k))_{k\geq 0}$ be a sequence that converges to $M\in \Omega_u$. Up to a subsequence, the sequence $(O(t_k))_{k\geq 0}$ converges to some $O\in SO_d$. The polar decomposition of $M$ can therefore be chosen equal to $OS_u$. In what follows we define by $\mathcal{O}_u$ the $\omega$-limit set of $\{O(t);\ t\geq 0\}$, and it is clear that
$$
\Omega_u=\{OS_u;\ O\in \mathcal{O}_u\}=\mathcal{O}_uS_u
$$

We can therefore state:

\newtheorem{squareroot}{Theorem}
\begin{squareroot}\label{squareroot}
Let $\Phi_u(t)=O(t)S(t)$ be the polar decomposition of $\Phi_u(t)$. The function
$$
t\longmapsto\Phi_u(t)^T\Phi_u(t)=S^2(t)
$$
converges when $t\longmapsto +\infty$ to the limit $S_u^2$ where $S_u$ is a symmetric nonnegative matrix, and is the common value of the square roots of $M^TM$ for $M\in \Omega_u$.

The $\omega$-limit set of $\{O(t);\ t\geq 0\}$ is a compact and connected subset of $SO_d$ denoted by $\mathcal{O}_u$ and $\displaystyle \Omega_u=\{OS_u;\ O\in \mathcal{O}_u\}=\mathcal{O}_uS_u$.

Moreover the switched system is asymptotically stable for the input under consideration if and only if $S_u=0$.
\end{squareroot}
\demo Everything has been proved, except the last assertion, which is obvious.

\hfill $\Box$
\vskip 0.2cm

\noindent \textbf{Remark} For a single matrix $B$, the polar decomposition is $e^{tB}=O(t)S(t)$ and it is clear that $B$ is Hurwitz if and only if $S(t)$ tends to $0$ as $t\longmapsto +\infty$.

\vskip 0.2cm

The matrix $S_u$ being the limit of $S(t)$ as $t$ tends to $+\infty$, it can sometimes be computed using a suitable sequence $(t_k)_{k\geq 0}$ (see Example \ref{exampleSu}). However it is also useful to write $S_u$ as an integral on $[0,+\infty[$. The convergence of this integral is a key point for the proof of Theorem \ref{tjrsvrai} in Section \ref{fastinputs}.

For $t\in [a_n,a_{n+1}[$ we have
$$
\Phi_u(t)=\exp((t-a_n)B_{u_n})\Phi_u(a_n) \quad\mbox{ and }\quad
\frac{d}{dt}\Phi_u(t)=B_{u_n}\Phi_u(t).
$$
Therefore
$$
\frac{d}{dt}S^2(t)=\frac{d}{dt}\Phi_u(t)^T\Phi_u(t)=\Phi_u(t)^T (B_{u(t)}^T+B_{u(t)} )\Phi_u(t)
$$
and
$$
S^2(t)=I+\int_0^t \Phi_u(s)^T (B_{u(s)}^T+B_{u(s)} )\Phi_u(s)ds.
$$
Since $S(t)$ decreases to $S_u$, the integral is convergent and
\begin{equation}
S_u^2=I+\int_0^{+\infty} \Phi_u(s)^T (B_{u(s)}^T+B_{u(s)} )\Phi_u(s)ds. \label{theEquation}
\end{equation}
Notice that for any $x\in\R^d$ and for any $t\geq 0$ the real number $x^T\Phi_u(t)^T (B_{u(t)}^T+B_{u(t)} )\Phi_u(t)x$ is nonpositive. The integral
\begin{equation}
\int_0^{+\infty} x^T\Phi_u(s)^T (B_{u(s)}^T+B_{u(s)} )\Phi_u(s)xds \label{theEquation2}
\end{equation}
is therefore absolutely convergent.


\section{Convergence on intervals}\label{interval}

This section is devoted to a technical result of uniform convergence of the function $\Phi_u$ on some sequences of intervals. Let $(t_k)_{k\geq 0}$ be an increasing (up to $+\infty$) sequence of positive numbers and $s>0$. Consider the sequence $(\Phi_k)_{k\geq 0}$ of functions from $[0,s]$ into $\mathcal{M}(d,\R)$ defined by
$$
\Phi_k(t)=\Phi_u(t_k+t).
$$
We are going to prove that this sequence of functions satisfies the hypothesis of Ascoli's Theorem. First of all the domain of the functions $\Phi_k$'s is $[0,s]$, a compact set. On the second hand they are uniformly bounded. Indeed for $i=1,\dots,p$ and for all $x\in \R^d$, the function $\displaystyle t\longmapsto \norm{e^{tB_i}x}$ is non increasing and consequently $\displaystyle \norm{e^{tB_i}}\leq 1$. Since for $t\in [a_n,a_{n+1}[$
$$
\Phi_u(t)=\exp((t-a_n)B_{u_n})\exp(\delta_{n-1}B_{u_{n-1}})\dots\exp(\delta_1B_{u_1})\exp(\delta_0 B_{u_0})
$$
we have as well $\displaystyle \forall t\geq 0,\  \norm{\Phi_u(t)}\leq 1$. This proves that the $\Phi_k$'s are uniformly bounded by $1$.

To finish we have to show that the sequence $(\Phi_k)_{k\geq 0}$ is equicontinuous. As
$$
\frac{d}{dt}e^{tB_i}=B_ie^{tB_i}, \quad \mbox{we have } \quad \forall t\geq 0 \quad \norm{\frac{d}{dt}e^{tB_i}}\leq \norm{B_i},
$$
and by the mean value theorem
$$
\forall\ s,t\geq 0 \qquad \norm{e^{tB_i}-e^{sB_i}}\leq \abs{t-s}\norm{B_i}.
$$
Consequently the function $\Phi_u$ is $\lambda$-Lipschitzian with $\lambda=\max\{\norm{B_i};\ i=1,\dots,p\}$. The functions $\Phi_k$ are as well $\lambda$-Lipschitzian, and the family $(\Phi_k)_{k\geq 0}$ is equicontinuous.

We can state:

\newtheorem{ascoli}[P1]{Proposition}
\begin{ascoli} \label{ascoli}
For the sequence $(\Phi_k)_{k\geq 0}$ defined above there exists a subsequence that converges uniformly to a continuous function
$$
t\longmapsto \Psi(t)
$$
from $[0,s]$ into $\Omega_u$.
\end{ascoli}

This proposition will be used in Sections \ref{SectChattering} and \ref{pairs}.


\section{The sets $\Omega_u$ and $\omega_u$}

We define by $\omega_u$ the set of points of $\Omega_u$ which are limits of sequences of switching times:
$$
\omega_u=\{M\in \Omega_u; \exists (a_{n_k})_{k\geq 1} \mbox{ such that } M=\lim_{k\mapsto+\infty}\Phi_u(a_{n_k})\}
$$
Recall that for each $B_i$ the set $\V_i$ is defined by
$$
\V_i=\{x\in\R^d;\ \forall t\in \R\ \norm{e^{tB_i}x}=\norm{x}\},
$$
and let us state the relation between the sets $\Omega_u$, $\omega_u$ and the $\V_i$'s.

\newtheorem{petitomega}[squareroot]{Theorem}
\begin{petitomega} \label{petitomega}
For any $x\in \R^d$ we have
$$
\Omega_u x\subset \omega_u x\bigcup \left( \cup_{i=1}^p\V_i\right)
$$
\end{petitomega}

\demo

Let $l\in \Omega_ux$. There exists an increasing (up to $+\infty$) sequence $(t_k)_{k\geq 1}$ such that $(\Phi_u(t_k)x)_{k\geq 1}$ converges to $l$. Let $a_{n_k}$ be the unique switching time such that $t_k\in[ a_{n_k},a_{n_k+1}[$. Up to a subsequence we can assume that $(\Phi_u(a_{n_k})x)_{k\geq 1}$ is also convergent, and denote its limit by $l_0\in \omega_ux$. The set of matrices $B_i$ being finite we can moreover assume that $u_{n_k}$ is for all $k$ equal to the same index $i$ for some $i\in\{1,...,p\}$.

Let us consider the various possibilities:
\begin{enumerate}
	\item If some subsequence of $(t_k-a_{n_k})_{k\geq 1}$ converges to $0$ as $k\mapsto +\infty$, then the uniform continuity of $\Phi_u$ implies $l=l_0\in \omega_u x$.
	\item If no such subsequence converges to $0$ then there exist $\delta>0$ and an integer $k_0$ such that $\forall k\geq k_0$, $t_k-a_{n_k}\geq \delta$. In that case, and for $k\geq k_0$,
	$$
	\Phi_u(a_{n_k}+\delta)x=e^{\delta B_i}\Phi_u(a_{n_k})x\longmapsto_{k\mapsto +\infty}e^{\delta B_i}l_0.
	$$
	It follows that $e^{\delta B_i}l_0$ and $l_0$ are two $\omega$-limit points for $x$.
	Therefore the equality $\norm{e^{\delta B_i}l_0}=\norm{l_0}$ holds (Proposition \ref{P1}), and $l_0\in \V_i$ according to Lemma \ref{decomposition}. To finish
	$$\begin{array}{ll}
	l & =\lim_{k\mapsto +\infty} e^{(t_k-a_{n_k})B_i}\Phi_u(a_{n_k})x\\
	  & =\lim_{k\mapsto +\infty} e^{(t_k-a_{n_k})B_i}l_0+\lim_{k\mapsto +\infty} e^{(t_k-a_{n_k})B_i}(\Phi_u(a_{n_k})x-l_0).
	\end{array}
	$$
	The second limit vanishes because $\Phi_u(a_{n_k})x-l_0$ converges to $0$ and\\ $\norm{e^{(t_k-a_{n_k})B_i}}\leq 1$. Since $l_0\in \V_i$ and $\V_i$ is $B_i$-invariant, the point $e^{(t_k-a_{n_k})B_i}l_0$ belongs to $\V_i$ for $k\geq k_0$ and so does the first limit. This shows that $l\in \V_i$ and completes the proof.
\end{enumerate}

\hfill $\Box$

\vskip 0.2cm

\noindent{\bf Remark}.
Let $M=\lim_{k\mapsto+\infty}\Phi_u(a_{n_k})\in \omega_u$. As noticed in the proof there exists a subsequence $(a_{n_{k_j}})$ for which the switching signal is constant, equal to some $i\in\{1,...,p\}$. Thus $\omega_u$ is the set of $\omega$-limit points that are obtained by sequences of switching times for which the input takes a unique value.

\vskip 0.2cm

From the second item of the proof, we can deduce the following proposition.

\newtheorem{chepa}[P1]{Proposition}
\begin{chepa} \label{chepa}
Let $x\in \R^d$, and let $(a_{n_k})_{k\geq 0}$ be a sequence of switching times such that $\Phi_u(a_{n_k})x$ converges to $l\in\Omega_ux$. If
\begin{itemize}
	\item there exists $i\in\{1,...,p\}$ such that $\forall k\geq 0$, $u_{n_k}=i$,
	\item there exists $\delta>0$ such that $\forall k\geq 0$, $\delta_{n_k}\geq \delta$,
\end{itemize}
then the limit $l$ belongs to $\V_i$.
\end{chepa}


\section{Stability for Regular Switching Signals}\label{regss}

The purpose of this section is to state and prove asymptotic stability results for what we call regular inputs (see Section \ref{SlowS} for the definition). One of their properties is to satisfy for all $x\in \R^d$
$$
\Omega_u x\subset \cup_{i=1}^p\V_i.
$$
The first task is to show that the inputs which do not verify this property are "chaotic", in the sense defined in the next subsection.

\subsection{Chaotic inputs}\label{SectChattering}
Let us assume that the sequence $(\Phi(a_{n_k}))_{k\geq 0}$ is convergent, and that the limit $l$ of $\Phi(a_{n_k})x_0$ does not belong to
$\displaystyle \cup_{i=1}^p \V_i$ for some $x_0$.

Pick $\tau>0$ and define for $k\geq 0$
$$
\Psi_k(t)=\Phi(a_{n_k}+t) \qquad \mbox{for } \ t\in [-\tau,\tau].
$$
By virtue of Proposition \ref{ascoli}, and up to a subsequence, we can assume that the sequence $(\Psi_k)_{k\geq 0}$ converges uniformly to a continuous function $\Psi$ from $[-\tau,\tau]$ into $\Omega_u$.

Moreover $\displaystyle \cup_{i=1}^p \V_i$ is closed and $\Psi([-\tau,\tau])$ compact, hence we can also assume that
\begin{equation}\label{vide}
\Psi([-\tau,\tau])\bigcap\ \cup_{i=1}^p \V_i=\emptyset.
\end{equation}
Let us make the hypothesis that there exist a sequence $(a_{n_l})_{l\geq 0}$ and $\delta>0$ such that one of the following conditions hold:
$$
[a_{n_l},a_{n_l}+\delta_{n_l}]\subset [a_{n_k}-\tau,a_{n_k}+\tau] \quad\mbox{ and }\ \delta_{n_l}\geq \delta>0
$$
or
$$
a_{n_l}< a_{n_k}+\tau< a_{n_l}+\delta_{n_l}\quad\mbox{ and }\ a_{n_k}+\tau-a_{n_l}\geq \delta>0.
$$
Up to a subsequence we can assume that $u_{n_l}$ is constant equal to $i$. But in that case, and according to Proposition \ref{chepa}, the limit of the sequence $(\Phi_u(a_{n_l})x_0)_{l\geq 0}$ belongs to $\V_i$, in contradiction with (\ref{vide}).

This discussion motivates the following definition, and proves Proposition \ref{nochattering}.

\newtheorem{chattering}{Definition}
\begin{chattering}\label{chattering}
The input $u$ is said to be chaotic if there exists a sequence $[t_k,t_k+\tau]_{k\geq 0}$ of intervals that satisfies the following conditions
\begin{enumerate}
	\item $\displaystyle t_k\longmapsto_{k\mapsto +\infty}+\infty$ and $\tau>0$.
	\item For all $\epsilon >0$ there exists $k_0$ such that for all $k\geq k_0$, the input $u$ is constant on no subinterval of $[t_k,t_k+\tau]_{k\geq 0}$ of length greater than or equal to $\epsilon$.
\end{enumerate}
An input that does not satisfy these conditions is called a \textit{non chaotic input}.
\end{chattering}

\newtheorem{nochattering}[P1]{Proposition}
\begin{nochattering}\label{nochattering}
If the input $u$ is non chaotic then for all $OS_u\in\Omega_u$:
$$
\mbox{Im}(OS_u)\subset \bigcup_{i=1}^p \V_i.
$$
\end{nochattering}


\subsection{Regular Inputs}\label{SlowS}

\newtheorem{hachedei}[chattering]{Definition}
\begin{hachedei} \label{hachedei}
The input $u$ is said to satisfy the assumption $H(i)$ if there exist a subsequence $(a_{n_k})_{k\geq 0}$ and $\delta>0$ such that
$$
\forall k\geq 0 \qquad u_{n_k}=i\ \mbox{ and }\ \delta_{n_k}\geq \delta.
$$
\end{hachedei}

From Proposition \ref{chepa} we deduce at once:

\newtheorem{slowswitch}[P1]{Proposition}
\begin{slowswitch}\label{slowswitch}
If the input $u$ satisfies the assumption $H(i)$ then:
$$
\forall x\in\R^d \qquad \Omega_ux\bigcap \V_i\neq \emptyset.
$$
\end{slowswitch}

We can know define what we call a "regular input", or "regular switching signal":

\newtheorem{snow}[chattering]{Definition}
\begin{snow} \label{snow}
An input $u$ is said to be {\bf regular} if it is non chaotic and satisfies the assumption $H(i)$ for $i=1,\dots,p$.
\end{snow}

\noindent{\bf Remark}.
The simplest non chaotic inputs are those for which the durations $\delta_n$ have a minimum $\delta>0$ (they are often called "slow switching inputs" in the literature). For these switching signals the hypothesis $H(i)$ is the consequence of one of the following weaker assumptions
\begin{enumerate}
	\item $m\{t\geq 0\;\ u(t)=i\}=+\infty\qquad$ ($m$ stands for the Lebesgue measure on the real line);
	\item an infinite number of $u_n$ take the value $i$.
\end{enumerate}
The paper \cite{HM} considers switching signals with average dwell-time: there exist $N_0>$ and 
$\tau_a>$ such the number $N_u(T,T+t)$ of discontinuities of $u$ in the interval $[T,T+t]$ satisfies:
$$
N_u(T,T+t)\leq N_0+\frac{t}{\tau_a}.
$$
Such inputs are clearly non chaotic.


\subsection{The stability theorems for regular inputs}

Recall that $S_r$ stands for the sphere of radius $r$ in $\R^d$.

\newtheorem{main}[squareroot]{Theorem}
\begin{main} \label{main}
If the sets $\V_i$ satisfy the condition:

\begin{itemize}
	\item[(C)] for $r>0$, no connected component of the set
$\displaystyle (\bigcup_{i=1}^p \V_i)\bigcap S_r$ $\qquad\qquad$ intersects all the $\V_i$'s;
\end{itemize}
then for every regular input $u$ the matrix $S_u$ is equal to $0$, and the switched system is asymptotically stable.
\end{main}

\demo Let $u$ be a regular input, and let us assume that $S_u$ does not vanish. Then there exists $x\in \R^d$ for which
$$
\Omega_u x\subset S_r
$$
for some $r>0$. Consider the two following facts:
\begin{enumerate}
	\item As $u$ satisfies $H(i)$ for $i=1,\dots,p$ we know by Proposition \ref{slowswitch} that $\displaystyle \Omega_u x\bigcap \V_i\neq \emptyset$;
	\item As $u$ is a non chaotic, we know by Proposition \ref{nochattering} that $\displaystyle \Omega_u x\subset \bigcup_{i=1}^p \V_i$.
\end{enumerate}
But $\Omega_u x$ is connected, and according to the two mentionned facts, it is a connected subset of $\displaystyle (\bigcup_{i=1}^p \V_i)\bigcap S_r$ that intersects all the $\V_i$.
This is in contradiction with Condition (C) hence $S_u=0$.

\hfill $\Box$

The assumption of Theorem \ref{main} implies obviously that $\displaystyle \bigcap_{i=1}^p \V_i^N=\{0\}$. Whenever that last holds we can find nice conditions, that is easy to check conditions, that imply Condition (C).

\newtheorem{corollaire}[squareroot]{Theorem}
\begin{corollaire} \label{corollaire}
Under the hypothesis
$$
\displaystyle \bigcap_{i=1}^p \V_i=\{0\},
$$
the matrix $S_u$ is equal to $0$, and the switched system is asymptotically stable for every regular input as soon as one of the following conditions hold:
\begin{enumerate}
  \item there exists $i$ such that $\dim \V_i=0$
  \item there exists $i$ such that $\dim \V_i=1$, and $\V_i\nsubseteq \V_j$, $\forall j\neq i$
	\item $p=2$
	\item $p>2$, but $\dim(\sum_{i=1}^p \V_i)>\sum_{i=1}^p \dim(\V_i)-p+1$
\end{enumerate}
In particular in the plane, that is for $d=2$, at least one of these conditions is satisfied as soon as $\V_i\neq\R^2$ for $i=1,\dots,p$.
\end{corollaire}

\demo
We have only to show that the stated conditions imply Condition (C) as soon as $\displaystyle \bigcap_{i=1}^p \V_i=\{0\}$ holds. It is obvious except for the fourth one which can be proved by induction.

Let us assume that a connected component $W$ of $ (\bigcup_{i=1}^p \V_i)\bigcap S_r$ intersects all the $\V_i$. Let $i_1=1$. There exists an index $i_2\neq i_1$ such that $\displaystyle \V_{i_1}\cap \V_{i_2}\neq \{0\}$, hence such that
$$
\dim(\V_{i_1}+\V_{i_2})\leq \dim \V_{i_1}+\dim \V_{i_2}-1.
$$
Let us assume that for some $k$, $1\leq k\leq p-1$, and some indices $i_1,\dots, i_k$, the inequality
$$
\dim(\sum_{j=1}^k \V_{i_j})\leq \sum_{j=1}^k \dim(\V_{i_j})-k+1
$$
holds. Then there exists an index $i_{k+1}\neq i_1,\dots,i_k$ such that $\displaystyle (\sum_{j=1}^k \V_{i_j})\cap \V_{i_{k+1}}\neq \{0\}$, hence such that
$$
\dim(\sum_{j=1}^k \V_{i_j}+\V_{i_{k+1}})\leq \sum_{j=1}^k \dim(\V_{i_j})-k+1+\dim(\V_{i_{k+1}})-1
$$

\hfill $\Box$

\vskip 0.3cm

\noindent{\bf Remarks}
\begin{enumerate}
	\item Condition (C) is actually a projective condition and could be replaced by: \textit{either the projection of $ \bigcup_{i=1}^p \V_i$ into the projective space $\mathbb{P}(\R^d)$ is not connected, or one of the $\V_i$'s is equal to $\{0\}$}.
	\item The assumption $ \bigcap_{i=1}^p \V_i=\{0\}$ cannot be avoided. If not a point $x_0\neq 0$ in this intersection can be in the kernel of $B_i$ for all $i$, hence a fixed point for any input.
\end{enumerate}

In case where Condition (C) does not hold, we cannot conclude to asymptotic stability for regular switching inputs. However the following result holds without that condition.

\newtheorem{general}[P1]{Proposition}
\begin{general} \label{general}
\begin{enumerate}
	\item Let $u$ be a non chaotic input. Then
	$$
	\forall M\in \Omega_u \quad \exists i\in\{1,\dots,p\} \quad \mbox{Im}(M)\subset \V_i.
	$$
	\item Let $u$ be a regular input. Then
	$$
	\rank (S_u)\leq \min_{i=1,\dots,p}\dim \V_i
	$$
\end{enumerate}
\end{general}

\demo
\begin{enumerate}
	\item Let $u$ be a non chaotic input and $M\in \Omega_u$. We know by Proposition \ref{nochattering} that $\displaystyle \mbox{Im}(M)\subset \bigcup_{i=1}^p \V_i$. But $\mbox{Im}(M)$ is a subspace of $\R^n$ and is therefore included in one of the $\V_i$'s.
	\item If $u$ satisfies moreover $H(i)$ for all $i=1,\dots,p$, then by virtue of Proposition \ref{chepa}
	$$
	\forall i=1,\dots,p \quad \exists M\in\Omega_u \quad \mbox{such that} \quad \mbox{Im}(M)\subset \V_i.
	$$
	Indeed if $M$ is the limit of $(\Phi_u(a_{n_k}))_{k\geq 0}$, where $\forall k\geq 0$ $u_{n_k}=i$ and $\delta_{n_k}\geq \delta>0$, then $Mx\in \V_i$ for all $x\in\R^d$.
	But for each $M\in \Omega_u$ there exists $O\in\mathcal{O}_u$ such that $M=OS_u$, and $\rank(S_u)=\rank(M)$.
\end{enumerate}

\hfill $\Box$


\section{Stability for chaotic inputs}\label{fastinputs}

In this section we deal with general inputs, that is with inputs which are possibly chaotic. For this purpose, we introduce the subspaces of $\R^d$
$$
\K_i=\ker(B_i^T+B_i) \qquad\mbox{ for }\quad i=1,...,p
$$
The subspace $\K_i$ is exactly the set of points $x$ for which the derivative of $\displaystyle \norm{e^{tB_i}x}^2$ vanishes at $t=0$. According to Lemma \ref{decomposition} the subspace $\V_i$ of $\K_i$ is $B_i$-invariant, but $\K_i$ is not, unless $\V_i=\K_i$. Notice that $B_i^T+B_i$ is negative out of $\K_i$, and that given a compact set $K$ that does not intersect $\K_i$, there exists $a>0$ such that
$$
\forall x\in K\qquad \qquad x^T(B_i^T+B_i)x\leq -ax^Tx
$$
It is proved in \cite{SVR} that for all $x\in\R^d$, the $\omega$-limit set $\Omega_ux$ is included in the union of the $\K_i$'s, {\it under the condition that the input $u$ satisfies}:
$$
\forall i\in\{1,\dots,p\} \qquad m\{t\geq 0;\ u(t)=i\}=+\infty.
$$
(In that paper the matrices $B_i$ are assumed to vanish on $\V_i$ but the proof works in our more general case).

However this requirement appears unnecessary: the sets $\Omega_ux$ are always included in the union of the $\K_i$'s, and actually an even better result holds. Let $J_u$ be the subset of $\{1,\dots,p\}$ defined by
$$
i\in J_u\Longleftrightarrow m\{t\geq 0;\ u(t)=i\}=+\infty,
$$
and let
$$
F_u=\bigcup_{i\in J_u}\K_i \subset \bigcup_{i=1}^p\K_i.
$$
Of course $J_u$ and $F_u$ depend on the input, but $F_u$ is always included in the union of the $\K_i$'s.

\newtheorem{tjrsvrai}[squareroot]{Theorem}
\begin{tjrsvrai} \label{tjrsvrai}
For all $x\in\R^d$, the set $\Omega_u(x)$ is included in $F_u$ and moreover verifies
$$
\forall i\in J_u \qquad\qquad \Omega_u(x)\bigcap \K_i \neq \emptyset
$$
\end{tjrsvrai}

\demo 

Let us assume that for some sequence $(t_k)_{k\geq 0}$ and some $x\in \R^d$
$$
\Phi_u(t_k)x\longmapsto l
$$ 
but that $l$ does not belong to $F_u$. This set being closed we can find $\epsilon>0$ and $\alpha>0$ such that
$$
\forall y\in B(l,2\epsilon) \quad \forall i\in J_u \qquad y^T(B_i^T+B_i)y\leq -\alpha y^Ty.
$$
Moreover, there exists $\eta>0$ such that
$$
\Phi_u(t)x\in B(l,\epsilon) \Longrightarrow \Phi_u(t+s)x\in B(l,2\epsilon) \qquad \mbox{ for }\quad \abs{s}\leq \eta.
$$
It is easy to see that the derivative of $\Phi_u(t)x$ is bounded by $\lambda(\norm{l}+2\epsilon)$,
where $\lambda=\max\{\norm{B_i};\ i=1,\dots,p\}$ is the Lipschitz constant of $\Phi_u$ (see Section \ref{interval}), as long as $\Phi_u(t)x$ belongs to $B(l,2\epsilon)$. The constant $\eta$ can therefore be chosen equal to $\frac{\epsilon}{\lambda(\norm{l}+2\epsilon)}$.

Now there exists an integer $k_0$ such that $\Phi_u(t_k)x\in B(l,\epsilon)$ as soon as $k\geq k_0$. We can also assume that the intervals $[t_k-\eta,t_k+\eta]$ are mutually disjoint. Denoting by $\chi_{\{u(s)\in J_u\}}$ the indicatrix function of the set $\{u(s)\in J_u\}$, we get
$$
\begin{array}{ll}
 & \displaystyle\int_0^{+\infty} x^T\Phi_u(s)^T (B_{u(s)}^T+B_{u(s)} )\Phi_u(s)x\ ds\\
&\\
\leq & \displaystyle\sum_{k\geq k_0}\int_{t_k-\eta}^{t_k+\eta} x^T\Phi_u(s)^T (B_{u(s)}^T+B_{u(s)} )\Phi_u(s)x\ ds\\
&\\
\leq & \displaystyle\sum_{k\geq k_0}\int_{t_k-\eta}^{t_k+\eta}(-\alpha)x^T\Phi_u(s)^T \Phi_u(s)x \chi_{\{u(s)\in J_u\}}\ ds\\
&\\
\leq & \displaystyle-\alpha\sum_{k\geq k_0}\int_{t_k-\eta}^{t_k+\eta}x^TS^2(s)x \chi_{\{u(s)\in J_u\}}\ ds\quad \leq 0.
\end{array}
$$
According to Section \ref{lift} all these integrals are convergent. In the last one the integration is made on a set of infinite measure, and $x^TS^2(s)x$ decreases to $x^TS_u^2x$ as $s$ tends to infinity. It follows that the convergence is possible only if $x^TS_u^2x=0$. Finally $l=OS_ux$ for some $O\in\mathcal{O}_u$, and we obtain $l=0$ in contradiction with the hypothesis $l\notin F_u$.

This proves the first assertion of the theorem. The second one being established in \cite{SVR}, we give only a sketch of the proof: if the compact set $\Omega_ux$ does not intersect the closed set $\K_i$ (where $i\in J_u$), then the distance between those two sets is positive, and for $t$ large enough $\Phi_u(t)x$ stays away from $\K_i$. The same kind of considerations than above shows that the integral, computed on the set of infinite measure $\{u(s)=i\}$, is divergent, a contradiction. 

\hfill $\Box$

\vskip 0.2cm

Thanks to Theorem \ref{tjrsvrai}, we can state a geometrical result of stability for general inputs.

\newtheorem{fastswitch}[squareroot]{Theorem}
\begin{fastswitch} \label{fastswitch}
With the previous notations, it is assumed that
$$
\displaystyle \bigcap_{i\in J_u} \K_i=\{0\}.
$$
Then the matrix $S_u$ is equal to $0$, and the switched system is asymptotically stable, as soon as  for $r>0$ no connected component of
$$
(\cup_{i\in J_u} \K_i)\bigcap S_r
$$
intersects all the $\K_i$'s for $i\in J_u$. This condition is verified in the following cases
\begin{enumerate}
	\item the cardinal of $J_u$ is $2$;
	\item the cardinal $q$ of $J_u$ is larger than $2$,\\ but $\dim(\sum_{i\in J_u} \K_i)>\sum_{i\in J_u} \dim(\K_i)-q+1$;
	\item there exists $i\in J_u$ such that $\dim \K_i=0$;
	\item there exists $i\in J_u$ such that $\dim \K_i=1$, but $\K_i\nsubseteq \K_j$, $\forall j\in J_u,\ j\neq i$.
\end{enumerate}
In particular for $d=2$, at least one of these conditions is satisfied as soon as $\K_i\neq\R^2$ for all $i\in J_u$.
\end{fastswitch}

\demo Similar to the one of Theorem \ref{corollaire}.

\hfill $\Box$


\section{Stability of pairs of Hurwitz matrices}\label{pairs}

In this section we deal with the case where $B_1$ and $B_2$ are Hurwitz. For the seek of generality our result is stated for general Lyapunov matrices, not necessarily equal to the identity.

\newtheorem{pairhurwitz}[squareroot]{Theorem}
\begin{pairhurwitz} \label{pairhurwitz}
Let $B_1$ and $B_2$ be two $d\times d$ Hurwitz matrices, assumed to share a common, but not necessarily strict, Lyapunov matrix $P$.
Then the switched system is asymptotically stable for any input as soon as
$$
\K_1\bigcap \K_2=\{0\}
$$
where $\K_i=\ker\left((PB_iP^{-1})^T+PB_iP^{-1}\right)$ for $i=1,2$.
\end{pairhurwitz}

\demo

Let us first assume that $P=Id$, and let $u$ be an input. There are two possibilities:
\begin{enumerate}
	\item Either $m\{t\geq 0;\ u(t)=i\}=+\infty$ for $i=1$ and $i=2$. According to Theorem \ref{fastswitch} the switched system is asymptotically stable.
	\item Or the equality $m\{t\geq 0;\ u(t)=i\}=+\infty$ is verified for only one index and we can assume without loss of generality that $m\{t\geq 0;\ u(t)=2\}<+\infty$. Let $l$ be an $\omega$-limit point for $x$: $l=\lim_{k\mapsto +\infty}\Phi_u(t_k)x$. Up to a subsequence we can choose $\tau>0$ such that the sequence $\Phi_k$ defined by $\Phi_k(t)=\Phi_u(t_k+t)$ converges uniformly to a function $\Psi$ on $[0,\tau]$ (see Section \ref{interval}). As
	$$
	\dt \Phi_k(t)=B_{u(t_k+t)}\Phi_k(t)
	$$
	we have
	\begin{equation}
	 \Phi_k(t)=\Phi_k(0)+\int_0^t B_{u(t_k+s)}\Phi_k(s)ds  \label{integral}
	\end{equation}
	The hypothesis $m\{t\geq 0;\ u(t)=2\}<+\infty$ implies that $m\{t\in [t_k,t_k+\tau];\ u(t)=2\}$ tends to $0$ as $k$ tends to $+\infty$, and up to a subsequence we can assume that  $\displaystyle m\{t\in [t_k,t_k+\tau];\ u(t)=2\}\leq 2^{-k}$. Let us denote by $E$ the set
	$$
	E=\bigcap_{k\geq 0}\bigcup_{j\geq k}\{t\in [0,\tau];\ u(t_j+t)=2\}
	$$
	Clearly the Lebesgue measure of this set vanishes. Moreover if $t\in [0,\tau]\setminus E$, then $u(t_k+t)$ is equal to $1$ for $k$ large enough, and we conclude that $B_{u(t_k+t)}$ converges to $B_1$ for almost all $t\in [0,\tau]$. The Lebesgue Theorem can be applied to the integral (\ref{integral}), so that we get
	$$
	 \Psi(t)=\Psi(0)+\int_0^t B_1\Psi(s)ds \quad\mbox{ and }\quad  \Psi(t)=\exp(tB_1)\Psi(0).
	$$
This shows that $\Psi(t)l=\exp(tB_1)l$ belongs to $\Omega_ux$ for all $t\in [0,\tau]$. Therefore $l\in \V_1$ by the very definition of this set, but $B_1$ being Hurwitz, $\V_1=\{0\}$ and $l=0$. The switched system is therefore asymptotically stable.
\end{enumerate}
If $P$ is not the identity matrix, we can replace $B_1$ and $B_2$ by
$$
\widetilde{B}_i=PB_iP^{-1} \qquad \mbox{for } i=1,2.
$$
The identity matrix is clearly a common Lyapunov matrix for $\widetilde{B}_1$ and $\widetilde{B}_2$.

\hfill $\Box$


\section{Examples and applications}\label{ExandAp}
\subsection{The particular case $d=2$}

Consider a collection $B_1,\dots,B_p$ of $2\times 2$ matrices satisfying Condition (\ref{condition}). In order to apply Theorem \ref{corollaire} we make the additional hypothesis
\begin{enumerate}
	\item $\displaystyle \bigcap_{i=1}^p \V_i=\{0\}$
	\item $\displaystyle \dim(\V_i)\leq 1$ for $i=1,\dots,p$.
\end{enumerate}
The switched system is then asymptotically stable for all regular inputs. For each matrix $B_i$ there are two possibilities
\begin{enumerate}
	\item $\dim(\V_i)= 1$. In that case $B_i$ has distinct eigenvalues $0$ and $\alpha<0$. Consequently $\K_i=\V_i$.
	\item $\dim(\V_i)= 0$. In that case $B_i$ is Hurwitz and the dimension of $\K_i$ is $0$ or $1$.
\end{enumerate}
Therefore the $\K_i$'s satisfy $\displaystyle \dim(\K_i)\leq 1$ for $i=1,\dots,p$. In view of the first additional hypothesis they satisfy as well $\displaystyle \bigcap_{i=1}^p \K_i=\{0\}$, except in the very particular following case:
\begin{itemize}
	\item There is an index $i_0$ for which
\begin{enumerate}
	\item $\displaystyle \dim(\V_{i_0})=0$, $\displaystyle \dim(\K_{i_0})=1$
	\item and $\forall i=1,\dots,p$, $\K_i=\K_{i_0}$.
\end{enumerate}
\end{itemize}

Apart from this particular case the switched system is asymptotically stable for all well distributed inputs, i.e. for inputs that verify $m\{t\geq 0;\ u(t)=i\}=+\infty$ for $i=1,...,p$, even if they are chaotic.


\subsection{An example of computation of $S_u$}\label{exampleSu}

Consider the $3\times 3$ matrices
$$
B_1=\begin{pmatrix}-1&0&0\\0&0&-1\\0&1&0\end{pmatrix},\quad
B_2=\begin{pmatrix}-1&0&0\\0&-1&0\\0&0&0\end{pmatrix},\quad
B_3=\begin{pmatrix}-1&0&0\\0&0&0\\0&0&-1\end{pmatrix}
$$
We have in the natural coordinates $(x,y,z)$ of $\R^3$:
$$
\V_1=\{x=0\}, \qquad \V_2=\{x=y=0\}, \qquad \V_3=\{x=z=0\}
$$
so that
$$
\bigcap_{i=1}^3\V_i=\{0\} \quad \mbox{ and} \quad \bigcup_{i=1}^3\V_i=\{x=0\}.
$$
The exponentials of these matrices are
$$
e^{tB_1}=\begin{pmatrix}e^{-t}&0&0\\0&\cos t&-\sin t\\0&\sin t&\cos t\end{pmatrix},\quad
e^{tB_2}=\begin{pmatrix}e^{-t}&0&0\\0&e^{-t}&0\\0&0&1\end{pmatrix},\quad
e^{tB_3}=\begin{pmatrix}e^{-t}&0&0\\0&1&0\\0&0&e^{-t}\end{pmatrix}.
$$
Let $u=(a_n,u_n)_{n\geq 0}$ be the input defined by $\displaystyle a_n=n\frac{\pi}{2}$ and
$$
u_{4k}=u_{4k+2}=1 \qquad u_{4k+1}=2 \qquad u_{4k+3}=3.
$$
We want to compute the matrix $S_u$ for this input. We have only to choose a sequence $(t_k)_k$ for which the limit of $\Phi_u(t_k)$ exists and is easy to calculate.
A straightforward computation gives
$$
e^{\frac{\pi}{2}B_3}e^{\frac{\pi}{2}B_1}e^{\frac{\pi}{2}B_2}e^{\frac{\pi}{2}B_1}=
\begin{pmatrix}e^{-2\pi}&0&0\\0&-1&0\\0&0&e^{-\pi}\end{pmatrix}
$$
so that $\Phi_u(4k\pi)$ tends to $\displaystyle \begin{pmatrix}0&0&0\\0&1&0\\0&0&0\end{pmatrix}$ as $k$ tends to $+\infty$. The matrix $S_u$ is therefore equal to $\displaystyle \begin{pmatrix}0&0&0\\0&1&0\\0&0&0\end{pmatrix}$.



\begin{thebibliography}{99}
\newcommand{\auth}{}
\newcommand{\tit}{}
\newcommand{\jou}{}
\newcommand{\bff}{}

\bibitem{agr}
\auth{A.~A.~Agrachev and D.~Liberzon},
\tit{Lie-algebraic stability criteria for switched systems},
\jou{SIAM J. Control Optim.},
\bff{40} (2001), 253-269.

\bibitem{AL00} \auth{D.~Angeli and D.~Liberzon},
\tit{A note on uniform global asymptotic stability of nonlinear switched systems in triangular form},
\jou{in Proc. 14th MTNS} (2000).

\bibitem{bbm}
\auth{M. Balde, U. Boscain and P. Mason},
\tit{A note on stability conditions for planar switched systems},
\jou{International Journal of Control},
\bff{82: 10} (2009) 1882-1888.

\bibitem{sw-balde} \auth{ M Balde, U. Boscain}, \tit{Stability of planar
switched systems: the nondiagonalizable case}, Communications in
Pure and Applied Analysis, AIMS, Vol 7, (2008) 1-21.

\bibitem{sw-1}  \auth{U.~Boscain}, \tit{Stability of planar switched
systems: the
linear single input  case},
\jou{SIAM J. Control Optim.}, \bff{41} (2002), 89-112.

\bibitem{branicky98} \auth{M.S.~Branicky}, \tit{Multiple Lyapunov functions and other analysis tools for switched and hybrid systems},
\jou{IEEE Trans. Aut. Control} \bff{43} (1998) 31-45.

\bibitem{chsi}
\auth {Y. Chitour and M. Sigalotti,}
\tit {On the stabilization of persistently excited linear systems,}
\jou{arXiv 0810, 2122v3 , math.OC, 18 may 2009}.

\bibitem{hart}  \auth{P. Hartmann}, \tit{Ordinary Differential equations},
\jou { John Willey and  sons, (1964)}.

\bibitem{HM} \auth{J.P. Hespanha, A.S. Morse}, \tit{Stability of switched systems with average dwell-time}
\jou{in Proc. 38th IEEE Conf. on Decision and Control (1999) 2655-2660}.

\bibitem{liberzon-book}  \auth{D. Liberzon}, \tit{Switching in systems
and
control}, Volume in series  Systems \& Control: Foundations \& Applications,
Birkh\"auser, Boston, (2003).
\bibitem{Marle} \auth {C.M.~Marle}, \tit{Syst\`emes dynamiques-Une introduction}, Ellipses, 2003.
\bibitem{Morse93} \auth{A.S. Morse}, \tit{Dwell-time switching}, \jou{in Proc. 2nd European Control Conf. (1993) 176-181}.
\bibitem{SVR}
\auth{U.~Serres, J.C.~Vivalda, P.~Riedinger},
\tit{On the convergence of linear switched systems},
\jou{hal-00371104 v1},
(2009).
\bibitem{vra} \auth{I. Vrabie}, \tit{Differential Equations,
An Introduction to Basic Concepts, Results and Applications}
\jou{World Scientific Publishing Co} (2004).
\end{thebibliography}
\end{document}